\documentclass[12pt]{article}

\usepackage{amsmath,amssymb,amsthm,amsfonts,sectsty,enumerate,titlesec,pifont,tikz,mathrsfs}
\usetikzlibrary{matrix,arrows,decorations.pathmorphing}

\newcommand{\vs}{\vspace{1pc}\noindent}
\newcommand{\hs}{\hspace{.2pc}\noindent}

\newcommand{\nin}{\noindent}

\newcommand{\aaa}{\ensuremath{\mathfrak{A}}} 

\titleformat{\section}{\bf\normalsize}{\thesection.}{3pt}{}

\titleformat{\subsection}[runin]{\bfseries\normalsize}{\thesubsection.}{3pt}{}

\theoremstyle{plain} 
\newtheorem{thm}{\bf Theorem}

\theoremstyle{definition} 

\subsectionfont{\normalsize}

\begin{document}

\setcounter{equation}{0}
\setcounter{thm}{0}
\setcounter{section}{-1}
\pagestyle{plain}


\let\thefootnote\relax\footnote{$\dag$ Email: pmxdcmo@nottingham.ac.uk. Telephone: 447592886332}

\vspace*{2pc}
\vs
\Large
Essential Spectra of Induced Operators on
\\
Subspaces and Quotients.

\normalsize
\vs
D.C. Moore$^\dag$
\vs

\vs\small
\textbf{Abstract.} Let $X$ be a complex Banach space and let $T$ be a bounded linear operator on $X$. For any closed $T$-invariant subspace $F$ of $X$, $T$ induces operators $T_{|F}:F \longrightarrow F$ and $T/F:X/F\longrightarrow X/F$. In this note, we give a simple proof of the fact that the essential spectra of $T_{|F}$ and $T/F$ are always contained in the polynomial hull of the essential spectrum of $T$.

\normalsize

\vs
\section{Introduction}

\nin
Henceforth, $X$ is a complex Banach space, $L(X)$ is the Banach algebra of bounded linear operators on $X$ and $K(X)$ is the (closed) ideal in $L(X)$ consisting of all compact operators on $X$. The ideal consisting of all finite-rank elements of $L(X)$ will be denoted using the symbol $K'(X)$. We adopt the convention whereby the {\em essential spectrum} of an operator $T\in L(X)$ is the set

$$
\sigma_{\mathrm{e}}(T) = \{ \lambda\in \mathbb{C}: \mbox{$\lambda-T$ is not a Fredholm operator}\}.
$$

\vs 
The {\em essential spectral radius}, $r_\mathrm{e}(T)$ is then taken to be $\sup\{ |\lambda|: \lambda\in \sigma_\mathrm{e}(T)\}$ when $X$ is infinite dimensional and $0$ when $X$ is finite dimensional. When $X$ is infinite dimensional, Atkinson's theorem asserts that $\sigma_{\mathrm{e}}(T)$ is the spectrum of the coset $T+K(X)$ in the Calkin algebra $L(X)/K(X)$. In this case, $\sigma_{\mathrm{e}}(T)$ is therefore always a non-empty compact subset of $\mathbb{C}$.

Given a closed subspace $F$ of $X$, $i_F$ is the inclusion operator sending $F$ into $X$, $Q_F$ is the quotient operator sending $X$ onto $X/F$ and $L_F(X)$ is the closed subalgebra of $L(X)$ consisting of all $T\in L(X)$ for which $TF\subseteq F$. It is clear that 

$$
L_F(X) = \{ T\in L(X): Q_FTi_F=0\}.
$$

\vs
The latest algebra is the natural domain of two continuous unital algbera homomorphisms. The first, which sends $T\in L_F(X)$ to $T_{|F}\in L(F)$, will be denoted $\pi_F^r$. The second, which sends $T\in L_F(X)$ to $T/F\in L(X/F)$, will be denoted $\pi_F^q$. It is straightforward to check that $\pi_F^r$ and $\pi_F^q$ map compact operators to compact operators.

\pagebreak

\nin
 It will also be convenient to set 

$$
\begin{array}{llllllll}
K_F(X) &=& L_F(X)\cap K(X) & ,& K_F'(X) & = & L_F(X)\cap K'(X)
\\[.5pc]
\aaa_F(X) & = & L_F(X)/K_F(X) &,& \aaa_F'(X) & = & L_F(X)/K_F'(X)
\end{array}
$$

\vs
The algebra $\aaa_F'(X)$ will play only the most minor r\^{o}le in these proceedings, so we need not be particularly troubled by the fact that it doesn't have a `canonical' topology (except when $X$ is finite dimensional).

Given any complex unital algebra $A$ and any element $T\in A$, we will write 

$$
\sigma^A(T) =\{\lambda\in \mathbb{C}: \mbox{$\lambda-T$ has no inverse in $A$}\},
$$

\vs
for the spectrum of $T$, suppressing the superscript $A$ when the algebra with respect to which the spectrum is being taken is unambiguous. The polynomially convex hull of a compact subset $S$ of $\mathbb{C}$ will be denoted $\widehat{S}$.

\vs
In this paper, we are concerned with the relationships between the sets 

$$
\begin{array}{llllllll}
\sigma_\mathrm{e}(T) & , & \sigma_\mathrm{e}(T_{|F}) &, &  \sigma_\mathrm{e}(T/F) &,& \sigma^{\aaa_F(X)}(T+K_F(X))
\\\\
&&&&& \mbox{and} & \sigma^{\aaa_F'(X)}(T+K_F'(X)).
\end{array}
$$

\vs
There is a some precedent for investigating questions of this type. For example, it is known that:

\vspace{.2pc}
\begin{enumerate}[(a)]
\item If $(U,V,W)$ is a permutation of $(T,T_{|F},T/F)$ then $\sigma_\mathrm{e}(U)\subseteq \sigma_\mathrm{e}(V)\cup \sigma_\mathrm{e}(W)$;
\item $\sigma^{\aaa_F'(X)}(T+K_F'(X)) = \sigma_\mathrm{e}(T)\cup \sigma_\mathrm{e}(T_{|F})$.
\end{enumerate}

\vspace{.2pc}

\nin
for any Banach space $X$, any closed subspace $F$ of $X$, and any $T\in L_F(X)$. The paper, \cite{SLAV}, of Djordjevi\'{c} and Duggal seems to be the most comprehensive open-access reference in this direction, and builds on the results obtained by Barnes in \cite{Barnes}. Equation (b) is proved in the monograph, \cite{BMSW}, of Barnes, Murphy, Smyth and West. All three of the references mentioned here contain various improvements of (a) and (b) in the presence of additional hypotheses on $X$, $F$ and $T$. For example, it is easy to show that:

\vspace{.2pc}

\begin{enumerate}[(c)]
\item $\sigma_{e}(T) = \sigma_\mathrm{e}(T_{|F})\cup \sigma_\mathrm{e}(T/F)$ if $F$ admits a $T$-invariant complement in $X$.
\end{enumerate}

\vspace{.2pc}

\nin
It is also clear that $\sigma_\mathrm{e}(T) = \sigma_\mathrm{e}(T_{|F})$ if $F$ is of finite codimension in $X$. The following example illustrates that the inclusion $\sigma_\mathrm{e}(T_{|F})\subseteq \sigma_\mathrm{e}(T)$ is too much to hope for in the general case. I am indebted to J. R. Partington for drawing my attention to this.

\pagebreak

\nin
\textbf{Example.} Given any open $U\subseteq \mathbb{R}$, let $L^2U$ denote the closure in $L^2(\mathbb{R})$ of the subspace $\{u\in C_c(\mathbb{R}): \mbox{$u$ has compact support in $U$}\}$ (identifying $C_c(\mathbb{R})$ with its image under the obvious embedding into $L^2(\mathbb{R})$ in the tradition fashional). Let $T:L^2(\mathbb{R})\longrightarrow L^2(\mathbb{R})$ be the operator which satisfies

$$
(Tu)(y) = u(y-1)
$$

\vs
for each $y\in \mathbb{R}$ and $u\in C_c(\mathbb{R})$. Plainly,

$$
T (L^2 (0,+\infty))\subseteq L^2(1,+\infty) \subseteq L^2(0,+\infty),
$$

\vs
so $F = L^2(0,+\infty)$ is a closed $T$-invariant subspace of $L^2(\mathbb{R})$. It is also clear that $TF\subseteq L^2(1,+\infty)$ is orthogonal to $L^2(0,1)$, so $T_{|F}$ is not lower semi-Fredholm. Thus $0\in\sigma_\mathrm{e}(T_{|F})$. Since $T$ is an invertible isometry, we also have $\sigma(T) \subseteq \mathbb{T}$, so $\sigma_\mathrm{e}(T_{|F})$ is certainly not contained in $\sigma_\mathrm{e}(T)$. 

\vs
One noticable feature of this example is that $\sigma_\mathrm{e}(T_{|F})$ is still contained in $\widehat{\sigma_\mathrm{e}(T)}$. The purpose of this paper is to give an easy proof that this is actually always true.

\section{The Result: Statement and Discussion}

\nin
The precise form of the result we obtain is as follows.

\begin{thm}
Let $X$ be a complex Banach space and let $F$ be an arbitrary closed subspace of $X$. Let $T\in L_F(X)$. Then 

$$
\sigma_\mathrm{e}(T_{|F})\cup \sigma_\mathrm{e}(T/F)\subseteq \widehat{\sigma_\mathrm{e}(T)}
$$

\vs
In particular, $\max\{ r_\mathrm{e}(T_{|F}), r_\mathrm{e}(T/F)\} \leq r_\mathrm{e}(T)$.
\end{thm}

\nin
The proof we give is based on two simple observations, both of which occur as exercises in the book, \cite{AaA}, of Abramovich and Aliprantis. These are:

\begin{enumerate}[(i)]
\item $T\in L_F(X)$ and $z\in \mathbb{C}\setminus \widehat{\sigma(T)}$ $\Rightarrow$ $(z-T)^{-1}\in L_F(X)$;
\item $T\in L_F(X)$ $\Rightarrow$ $\sigma(T_{|F})\subseteq \widehat{\sigma(T)}$.
\end{enumerate}

\nin
The details appear here for the convenience of the reader. Let $L(F,X/F)$ be the Banach space of all bounded linear operators from $F$ into $X/F$ and define 

$$
G: \mathbb{C}\setminus \widehat{\sigma(T)}\longrightarrow L(F,X/F)
$$

\vs
by setting $G(z) = Q_F((z-T)^{-1})i_F$ for each $z\in \mathbb{C}\setminus \widehat{\sigma(T)}$.

\pagebreak

\nin
It is clear that $G$ is holomorphic in $\mathbb{C}\setminus \widehat{\sigma(T)}$ and that assertion (i) is equivalent to $G$ being identically $0$ on $\mathbb{C}\setminus \widehat{\sigma(T)}$. By the identity theorem for Banach-space-valued functions, it suffices to show that $G(z)=0$ for $|z|>r(T)$. Let $z\in \mathbb{C}$ with $|z|>r(T)$. Then, assuming that $T\in L_F(X)$,

$$
G(z) = Q_F \left( \lim_{n\to \infty} \sum_{m=0}^n \dfrac{T^m}{z^{m+1}} \right)i_F = \lim_{n\to \infty}\sum_{m=0}^n \dfrac{Q_FT^mi_F}{z^{m+1}}=0,
$$

\vs
since $T^m$ belongs to $L_F(X)$ for each $m$. This clearly establishes that 

$$
\sigma^{L_F(X)}(T)\subseteq \widehat{\sigma(T)}
$$

\nin
Assertion (ii) and the fact that we also have $\sigma(T/F)\subseteq \widehat{\sigma(T)}$ follows immediately, since 

$$\begin{array}{lrlllllllllll}
&\sigma(T_{|F}) &=& \sigma(\pi_F^r(T)) &\subseteq & \sigma^{L_F(X)}(T), 
\\[.5pc]
\mbox{and}&\sigma(T/F) &=&\sigma(\pi_F^q(T))&\subseteq& \sigma^{L_F(X)}(T).
\end{array}$$

\vs
Observation (i) has another interesting consequence (which is not noted in \cite{AaA}). Let $\lambda$ be an isolated point of $\sigma(T)$ and suppose that there is a neighbourhood $V$ of $\lambda$ such that $V\setminus \{\lambda\} \subseteq \mathbb{C}\setminus \widehat{\sigma(T)}$. By choosing $r>0$ sufficiently small, we can arrange it that $\gamma(t) = \lambda + r\mathrm{e}^{it}$ belongs to $V$ for every $t\in [0,2\pi)$. The spectral projection 

\begin{equation}
\label{specproj}
P_T(\lambda) = \dfrac{1}{2\pi i}\int_\gamma (z-T)^{-1} \hs \mathrm{d}z
\end{equation}

\vs
is thus expressed as a limit of sums of elements $(z-T)^{-1}$ for $z\in \mathbb{C}\setminus \widehat{\sigma(T)}$. That $P_T(\lambda)\in L_F(X)$ follows. By (ii), either $\lambda-T_{|F}$ is invertible or $\lambda$ occurs as an isolated point of $\sigma(T_{|F})$. It therefore makes sense to consider the spectral projection of $T_{|F}$ associated with $\lambda$. Since $\pi_F^r$ is an algebra homomorphism, we have

$$
P_{T_{|F}}(\lambda) = \dfrac{1}{2\pi i}\int_\gamma (z-T_{|F})^{-1}\hs \mathrm{d}z = \dfrac{1}{2\pi i}\int_\gamma \pi_F^r((z-T)^{-1})\hs \mathrm{d}z
$$

\vs
As $\pi^r_F$ is a continuous linear map, this gives

$$
P_{T_{|F}}(\lambda) = \pi_F^r\left( \dfrac{1}{2\pi i}\int_\gamma (z-T)^{-1} \hs \mathrm{d}z \right) = \pi_F^r(P_T(\lambda))=P_T(\lambda)_{|F},  
$$

\vs
with an entirely parallel conclusion available for $P_{T/F}(\lambda)$.

\pagebreak

\nin
The proof of Theorem 1 also uses the following standard fact (c.f. \cite{Davies}, Lemma 4.3.17): if $\lambda\in \mathbb{C}$ is an isolated point of $\sigma(T)$ then $P_T(\lambda)$ has finite rank if and only if $\lambda-T$ is a Fredholm operator.

\section{Proof of Theorem 1.}

\nin
Let $\pi$ be either of $\pi_F^r$ or $\pi_F^q$, so that $\pi(T) = T_{|F}$ or $\pi(T)=T/F$. It is clear that $\pi$ is a continuous unital algebra homomorphism and maps finite rank operators to finite rank operators. Let $U=\mathbb{C}\setminus\widehat{\sigma_\mathrm{e}(T)}$. By the punctured neighbourhood theorem (\cite{Muller}, Theorem 19.4), $U\cap \sigma(T)$ consists of isolated points of $\sigma(T)$. Thus every $\lambda\in U$ has a neighbourhood, $V$, such that 

$$
V\setminus \{\lambda\} \subseteq \mathbb{C}\setminus \widehat{\sigma(T)}.
$$

\vs
Let $\lambda\in U$ and let $V$ be a neighbourhood of $\lambda$ with the property above. Since $\lambda\notin \sigma_\mathrm{e}(T)$, $P_T(\lambda)$ is a finite rank operator. Also, as we saw in Section 1, either $\lambda\in \mathbb{C}\setminus \sigma(T)$ or is an isolated point of $\sigma(\pi(T))$. In the first case, $\lambda-\pi(T)$ is certainly Fredholm and the proof ends here. Otherwise, the argument at the end of Section 1 gives $P_T(\lambda)\in L_F(X)$ and $P_{\pi(T)}(\lambda) = \pi(P_T(\lambda))$. Since $\pi$ maps finite rank operators, $\lambda-\pi(T)$ is Fredholm, so $\lambda\notin \sigma_\mathrm{e}(\pi(T))$. Having shown that $\lambda\in \mathbb{C}\setminus \widehat{\sigma_\mathrm{e}(T)}\Rightarrow \mathbb{C}\setminus \sigma_\mathrm{e}(\pi(T))$, when $\pi$ is either of $\pi_F^r$ or $\pi_F^q$, we conclude that 

$$
\sigma_\mathrm{e}(T/F)\cup \sigma_\mathrm{e}(T_{|F})\subseteq \widehat{\sigma_\mathrm{e}(T)},
$$

\vs
as asserted. The inequality $\max\{ r_\mathrm{e}(T_{|F}), r_\mathrm{e}(T/F)\} \leq r_\mathrm{e}(T)$ follows.

\section{Closing Remarks}

\nin
All of the ingredients of the proof have been known in isolation for many years. In many ways, it would be quite surprising if no version of it had appeared before. Another, similarly elementary result which does not appear to exist in the literature is the fact that, with $X,F$ and $T$ as in Theorem 1, 

$$
\sigma_\mathrm{e}(T_{|F})\cup \sigma_\mathrm{e}(T/F)\cup \sigma_\mathrm{e}(T)\subseteq \sigma^{\mathfrak{A}_F(X)}(T+K_F(X)).
$$

\vs
This `soft' result simply follows from the fact that $\pi_F^r$, $\pi_F^q$ and the inclusion of $L_F(X)$ into $L(X)$ map compact operators to compact operators. It is also not difficult to show, in the special case where $F$ is complemented in $X$, that $\sigma_\mathrm{le}(T_{|F})\cup \sigma_\mathrm{re}(T/F)\subseteq \sigma_\mathrm{e}(T)$, where $\sigma_\mathrm{le}$ and $\sigma_\mathrm{re}$ denote left and right essential spectra (definitions can be found in \cite{Muller}, page 172). This latest observation also appears to be missing from the literature.

\pagebreak

\nin
{\bf Acknowledgements:} I am indebted to the EPSRC (Grant EP/K503101/1) for their financial support during the preparation of this note. Some of the material in this paper may appear the author's PhD thesis. 

\vs

%
%
%
%
%
%
%
%
%
%
%
%
%

\end{document}